\documentclass[9pt,shortpaper,twoside,web]{ieeecolor}
\usepackage{generic}

\usepackage{hyperref}%
\hypersetup{colorlinks=true, linkcolor=blue, breaklinks=true, urlcolor=blue}
\usepackage{doi}

\usepackage[all]{xy}
\usepackage{amsmath}
\usepackage{amscd}
\usepackage{mathrsfs}
\usepackage{amssymb}
\usepackage{tikz}
\usepackage[utf8]{inputenc}
\usepackage[yyyymmdd]{datetime}
\usepackage{graphicx}
\usepackage{subcaption}
\usepackage{lipsum}
\usepackage{multicol} \usepackage{mathtools}
\usepackage{multirow}
\usepackage{threeparttable}
\usepackage{pdfpages}
\usepackage{stmaryrd}
\usepackage{selectp}

\usepackage{multicol} \usepackage{mathtools}
\usepackage{version,xspace}
\usepackage{url,doi}

 \newcommand{\setdef}[2]{\{#1
	\; | \; #2\}}

\newcommand\oprocendsymbol{\hbox{$\triangle$}}
\newcommand\oprocend{\relax\ifmmode\else\unskip\hfill\fi\oprocendsymbol}

\let\labelindent\relax
                  
\usepackage{enumitem} \renewcommand{\theenumi}{(\roman{enumi})}

\DeclareSymbolFont{bbold}{U}{bbold}{m}{n}
\DeclareSymbolFontAlphabet{\mathbbold}{bbold}
\newcommand{\vect}[1]{\mathbbold{#1}}

\newcommand{\real}{\mathbb{R}}

\newcommand{\seminorm}[1]{{\left\vert\kern-0.25ex\left\vert\kern-0.25ex\left\vert #1
		\right\vert\kern-0.25ex\right\vert\kern-0.25ex\right\vert}}

\newcommand{\semimeasure}[1]{\mu_{\seminorm{\cdot}}\kern-0.5ex\left(#1\right)}

\DeclareMathOperator{\Ker}{\mathrm{Ker}}

\renewcommand{\top}{\mathsf{T}} 

\newtheorem{theorem}{Theorem}[section]
\newtheorem{proposition}[theorem]{Proposition}

\newtheorem{lemma}[theorem]{Lemma}
\newtheorem{definition}[theorem]{Definition}
	\newtheorem{remark}[theorem]{Remark}
	\newtheorem{example}[theorem]{Example}

\newcommand{\suchthat}{\;\ifnum\currentgrouptype=16 \middle\fi|\;}
\newcommand{\scirc}{\raise1pt\hbox{$\,\scriptstyle\circ\,$}}

\def\BibTeX{{\rm B\kern-.05em{\sc i\kern-.025em b}\kern-.08em
    T\kern-.1667em\lower.7ex\hbox{E}\kern-.125emX}}
\begin{document}
\title{Monotonicity and Contraction on Polyhedral Cones}
\author{Saber Jafarpour, \IEEEmembership{Member, IEEE} and Samuel Coogan
  \IEEEmembership{Senior Member, IEEE}
\thanks{This work was supported in part by the National Science Foundation under grant 2219755 and the Air Force Office of Scientific Research under Grant FA9550-23-1-0303.}
\thanks{Saber Jafarpour is with the Department of Electrical and Computer Engineering, University of Colorado Boulder, USA {\tt\small saber.jafarpour@colorado.edu}}
\thanks{Samuel Coogan is with the School of Electrical and Computer Engineering, Georgia Institute of Technology, USA, {\tt\small sam.coogan@gatech.edu}}
}

\maketitle

\begin{abstract}
 In this note, we study monotone dynamical systems with respect to polyhedral cones. 
Using the half-space representation and the vertex representation, we propose three equivalent conditions to certify monotonicity of a dynamical system with respect to a polyhedral cone. 
We then introduce the notion of gauge norm associated with a cone and provide closed-from formulas for computing gauge norms associated with polyhedral cones.
A key feature of gauge norms is that contractivity of monotone systems with respect to them can be efficiently characterized using simple inequalities. This result generalizes the well-known criteria for Hurwitzness of Metzler matrices and provides a scalable approach to search for Lyapunov functions of monotone systems with respect to polyhedral cones. 
Finally, we study the applications of our results in transient stability of dynamic flow networks and in scalable control design with safety guarantees.  

\end{abstract}

\section{Introduction}

\paragraph*{Motivation and Problem Statement}


Monotone systems are a class of dynamical systems characterized by preserving a partial ordering along their trajectories. 
The framework of monotone systems has been successfully used to model complex systems in nature such as biochemical cascade reactions~\cite{EDS:07} as well as engineered system such as transportation networks~\cite{GC-EL-KS:15}. 
It is known that monotone systems exhibit highly ordered dynamical behaviors~\cite{HLS:95} that can be used to establish stability of their interconnection~\cite{DA-EDS:03}, to develop computationally efficient techniques for their control synthesis~\cite{AR:15,YE-DP-DA:17} and to perform reachability analysis to ensure their safety~\cite{SC-MA:15b}. 

 Contraction theory is a classical framework for studying dynamical systems where stability is defined incrementally between two arbitrary trajectories. Contracting systems feature desirable transient and asymptotic behaviors including i) forgetting their initial conditions, ii) exponential convergence to a single trajectory, and iii) input-to-state robustness with respect to disturbances and unmodelled dynamics. While the study of contracting systems can be traced back to the 1950s, many recent works have  focused on infinitesimal frameworks~\cite{WL-JJES:98} and Finsler-Lyapunov frameworks~\cite{FF-RS:14} for analysis of contracting systems.


A large body of the research in monotone system theory focuses on 
cooperative systems, i.e., systems that are monotone with respect to the positive orthant. 
It is well known that cooperative systems are amenable to efficient stability analysis using {\color{black}suitable} Lyapunov functions~\cite{GD-HI-AR-BSR:15}, a feature that can be used to develop computationally efficient techniques for control design of large-scale cooperative systems. For cooperative dynamical systems over networks, several recent works have focused on construction of sum-separable or max-separable Lyapunov functions, i.e., Lyapunov functions that can be represented as sum or maximum of scalar-valued functions of each state component~\cite{HF-BB-MJ:18,YK-BB-MC:20}.
It turns out that, for cooperative systems, contractivity plays an essential role in the design of sum-separable and max-separable Lyapunov functions. In~\cite{SC:19} contraction with respect to the $\ell_{\infty}$- and $\ell_{1}$-norm has been used to establish existence of sum-separable and max-separable Lyapunov function for cooperative systems. In~\cite{IRM-JJES:17} contraction theory with respect to a Riemannian metric has been used to study sum-separable and max-separable Lyapunov functions for cooperative systems.

Monotonicity with respect to arbitrary cones accommodates a significantly broader class of systems than cooperativity. For example, for linear dynamical systems, it is known that if all the eigenvalues of the system are real, then there exists a cone with respect to which the system is monotone~\cite[Theorem 3.5]{AB-RJP:94}. Furthermore, monotonicity with respect to polyhedral cones is closely-related to the problem of existence of an internally positive realization~\cite{YO-HM-SK:84,LF:96}. For linear systems, monotonicity with respect to arbitrary cones on $\real^n$ has been studies in~\cite{HS-MV:70}. The paper~\cite{FB-MF-EH:74} studies necessary and sufficient conditions for positivity of a linear operator with respect to a polyhedral cone. For nonlinear systems, a variational characterization of monotonicity with respect to arbitrary cones is presented in~\cite{DA-ES:08}, and a generalization of positivity with respect to an arbitrary cone is proposed in~\cite{FF-RS:16}.   

Nonetheless, certifying monotonicity of systems with respect to arbitrary cones is usually computationally complicated. Moreover, many techniques developed for cooperative systems, including those in \cite{GD-HI-AR-BSR:15, HF-BB-MJ:18,YK-BB-MC:20, SC:19, IRM-JJES:17}, do not generalize easily or at all to the broader class of monotone systems. In particular, the connection between contraction theory and monotone system theory with respect to arbitrary cones and the existence of {\color{black}suitable} Lyapunov functions for stability analysis of monotone systems with respect to arbitrary cones is not well understood or studied. Exceptions are~\cite{DK-FF:20}, which considers searching for a polyhedral cone which makes a nonlinear system monotone, and~\cite{YK-FF:22}, which studies incremental stability of monotone systems with respect to arbitrary or polyhedral cones.



\paragraph*{Contributions}
In this note, we study monotonicity and contractivity of dynamical systems with respect to polyhedral cones.
First, given an arbitrary polyhedral cone, we provide three equivalent characterization of dynamical systems that are monotone with respect to this polyhedral cone. To the best of our knowledge, this result is novel for non-proper or non-pointed polyhedral cones.
Second, given a proper cone and a vector in its interior, we introduce the notion of the gauge and the dual gauge norms as natural metrics for studying contractivity of monotone systems. We provide closed-form formulas for computing the gauge and dual gauge norms and characterize the sublevel sets of the gauge and dual gauge matrix measures.
Third, for monotone systems with respect to proper polyhedral cones, we provide necessary and sufficient condition for their contractivity with respect to both the gauge norm and the dual gauge norm. Our conditions for contractivity with respect to the gauge and dual gauge norms are generalizations of the closed-form expressions for $\ell_{\infty}$-norm and $\ell_{1}$-norm contractivity of cooperative systems.
 As the first application of our analysis, we study transient behavior of edge flows in networks with nodal dynamics. The key element of our approach is a non-pointed polyhedral cone defined using the structure of the network. Finally, using our results in this note, we propose a scalable control design scheme for a nonlinear dynamical system constrained to be in a safe subset of its state space. Our approach uses linear programming to design a suitable feedback controller such that the closed-loop system avoids the unsafe region.
 

\section{Notations and Mathematical Preliminary}

Let $\mathcal{L}$ be a set with a relation $\preceq$. Then $\preceq$ is a preorder if
\begin{enumerate}
    \item\label{p1:self} $x \preceq x$, for every $x\in \mathcal{L}$; 
    \item\label{p3:trans} $x\preceq y$ and $y\preceq z$ implies that $x\preceq z$.
\end{enumerate}
A preorder $\preceq$ is a partial order if it additionally satisfies
\begin{enumerate}\setcounter{enumi}{2}
   \item\label{p2:sym} $x \preceq y$ and $y \preceq x$ implies that $x=y$;
\end{enumerate}
Let $W$ be a vector spaces and $\mathcal{S}\subseteq W$. The convex hull of $\mathcal{S}$ is denoted by $\mathrm{conv}(\mathcal{S})$ and the interior of $\mathcal{S}$ is denoted by $\mathrm{int}(\mathcal{S})$. 
The standard partial ordering $\le $ on $\real^n$ is defined as $v\le w$ if $v_i\le w_i$, for every $i\in \{1,\ldots,n\}$. A matrix $A\in \real^{n\times n}$ is Metzler if all its off-diagonal elements are non-negative. The set of $n\times n$ Metzler matrices are denoted by $\mathbb{M}_n\subset \real^{n\times n}$. For a given matrix $A\in \real^{n\times m}$, we denote its transpose by $A^{\top}\in \real^{m\times n}$ and its Moore–Penrose inverse by $A^{\dagger}\in \real^{m\times n}$. Given a vector $\eta\in \real^n$, we define the diagonal matrix $\mathrm{diag}(\eta)\in \real^{n\times n}$ by $\mathrm{diag}(\eta)_{ii}=\eta_i$, for every $i\in\{1,\ldots,n\}$. For a vector space $W$, we denote its dual with $W^*$ and the dual-pairing between $W$ and $W^*$ is denoted by $\langle\phi, v\rangle = \phi(v)$, for every $\phi\in W^*$ and every $v\in W$. Given a norm $\|\cdot\|$ on the vector space $W$, the induced norm on the dual space $W^*$ is defined by $\|\phi\| = \max\setdef{|\langle\phi, v\rangle|}{\|v\|\le 1}$ and {\color{black}the \emph{dual norm} on $W$ is defined by $\|v\|^{\mathrm{d}} =\max\setdef{|\langle\phi, v\rangle|}{\|\phi\|\le 1}$.} A set $S$ is absorbent in the vector space $W$, if for every $v\in W$, there exists $r>0$ such that $c v\in S$, for every $c$ such that $|c|\le r$~\cite[Definition 4.1.2]{LN-EB:85}.
Given a seminorm $\seminorm{\cdot}$ on $\real^n$, its kernel is defined by $\Ker\seminorm{\cdot}=\setdef{x\in \real^n}{\seminorm{v}=0}$. The induced seminorm of $A$ is defined by $\seminorm{A} = \sup\setdef{\seminorm{Ax}}{\seminorm{x}=1, x \perp \Ker\seminorm{\cdot}}$ and the matrix semi-measure of $A$ with respect to the seminorm $\seminorm{\cdot}$ is defined by $\mu_{\seminorm{\cdot}}(A) = \lim_{h\to 0^+}\frac{\seminorm{I + hA}-1}{h}$~\cite[Definition 4]{SJ-PCV-FB:22}.
A control system $\dot{x}=f(x,u)$ is incrementally exponentially stable with rate $c\in \real_{\ge 0}$ if, for every norm $\|\cdot\|$, there exists $M\in \real_{\ge 0}$ such that
\begin{align*}
    \|x_u(t)-y_u(t)\|\le M e^{-ct} \|x_u(0)-y_u(0)\|,\quad\mbox{ for all }t\in \real_{\ge 0},
\end{align*}
for every two trajectories of the system $t\mapsto x_u(t)$ and $t\mapsto y_u(t)$ with the same input $t\mapsto u(t)$. The system is called semi-contracting with respect to a seminorm $\seminorm{\cdot}$ with rate $c\in \real_{\ge 0}$ if
\begin{align*}
    \seminorm{x_u(t)-y_u(t)} \le e^{-ct} \seminorm{x_u(0)-y_u(0)},\quad\mbox{ for all }t\in \real_{\ge 0},
\end{align*}
for every two trajectories of the system $t\mapsto x_u(t)$ and $t\mapsto y_u(t)$ with the same input $t\mapsto u(t)$.

\section{Cones, positive operators, and Metzler operators}

A non-empty subset $K\subseteq \real^n$ is a cone if (i) $|\lambda| K \subseteq K$, for every $\lambda\in \real$, (ii) $K$ is closed in $\real^n$, and (iii) $K$ is convex, i.e., $K+K\subseteq K$. A cone $K\subseteq \real^n$ is called pointed if $K\cap (-K) = \{0\}$ and is called proper if $\mathrm{int}(K)\ne \emptyset$. Given a cone $K\subseteq \real^n$, the preorder $\preceq_{K}$ on $\real^n$ is given by
       \begin{align*}
     x\preceq_K y \iff y-x\in K.
       \end{align*}
       If $K\subseteq \real^n$ is a pointed cone, then the preorder $\preceq_{K}$ is a partial order. For every $x\preceq_K y$, the interval $[x,y]_{K}$ is defined
       by $[x,y]_K=\setdef{z\in V}{x\preceq_K z \preceq_K y}$
       For $S\subseteq \real^n$, the polar set $S^*$ is 
    \begin{align*}
      S^* = \setdef{\phi\in (\real^n)^*}{ \langle \phi , x \rangle\ge 0, \mbox{ for all }x\in
      S}. 
      \end{align*} 
     The polar set $K^*$ is again a cone and is usually denoted by the dual cone of $K$. {\color{black}Moreover, if the cone $K$ is proper then the dual cone $K^*$ is pointed and if the cone $K$ is pointed then the dual cone $K^*$ is proper.}
A cone $K\subseteq \real^n$ is polyhedral if 
\begin{align}\label{eq:polyhedral}
K = \setdef{x\in \real^n}{\langle \phi_i, x \rangle \ge 0, \;\;\; \forall \;\; i\in \{1,\ldots,m\}}
\end{align}
where $\phi_i:\real^n\to \real$ is a linear functional for every $i\in\{1,\ldots,m\}$. Given a polyhedral cone $K\subseteq \real^n$, there exist two matrices $H\in \real^{m\times n}$ and $V\in \real^{n\times q}$ such that $K$ has the following equivalent representations~\cite[Theorem 1.3]{GMZ:12}: 
\begin{align}
    K &=\setdef{x\in \real^n}{Hx \ge \vect{0}_m},\label{eq:H}\\
    K &=\setdef{Vx\in \real^n}{x \ge \vect{0}_q}.\label{eq:V}
\end{align}
The representation~\eqref{eq:H} is called a \emph{half-space representation} ($H$-rep) of the cone $K$ and the matrix $H$ is called the \emph{representation matrix} for the cone $K$. The representation~\eqref{eq:V} is called a \emph{vertex representation} ($V$-rep) of the cone $K$ and the matrix $V$ is called the \emph{generating matrix} for the cone $K$~\cite{GMZ:12}. 
Given a polyhedral cone $K\in \real^n$ with an $H$-rep as in equation~\eqref{eq:H} and a $V$-rep as in equation~\eqref{eq:V}, the pair $(H,V)$ is called a representation for the polyhedral cone $K$. For a polyhedral cone $K$ with a representation $(H,V)$, the dual cone $K^*$ is a polyhedral cone with representation $(V^{\top},H^{\top})$~\cite[Theorem 4.7]{KF:20}. 
Several algorithms exists for transforming $H$-rep to $V$-rep and vice versa, including Fourier-Motzkin elimination~\cite{GMZ:12} and the double description method~\cite{KF-AP:95}. 
  
The next lemma uses a representation $(H,V)$ of the polyhedral cone $K$ to find the closed-form expressions for the preorders $\preceq_{K}$ and $\preceq_{K^*}$. We refer to~\cite{SJ-SC:22} for the proof of this lemma. 

\begin{lemma}[$H$-rep and $V$-rep of polyhedral cones]\label{thm:inequality}
Let $K\subset \real^n$ be a polyhedral cone with representation $(H,V)$ such that $H\in \real^{m\times n}$ and $V\in \real^{n\times q}$. The following statements are equivalent:
    \begin{enumerate}
    \item\label{p1:inequality} $v \preceq_K w$;
    \smallskip
     \item\label{p2:inequality} $Hv \le Hw$.
      \end{enumerate}
      Additionally, the following statements are equivalent:
      \begin{enumerate}\setcounter{enumi}{2}
    \item\label{p3:inequality} $v \preceq_{K^*} w$;
    \smallskip
     \item\label{p4:inequality} $V^{\top}v \le V^{\top}w$.
      \end{enumerate}
    \end{lemma}
\smallskip

    Given a cone $K\subseteq \real^n$, the linear operator $A:\real^n\to \real^n$ is called
           \begin{enumerate}
           \item\label{c1:positive} $K$-positive if $AK\subseteq K$; 
            \item\label{c2:metzler} $K$-Metzler if, for every $\phi\in K^*$ and every $v\in K$ such that
              $\langle \phi, v \rangle=0$, we have $\langle \phi, Av \rangle\ge 0$. 
\end{enumerate}
In the literature, $K$-Metzler matrices are sometimes referred to as \emph{cross-positive} matrices~\cite{HS-MV:70} or \emph{$K$-cooperative} matrices~\cite{MWH-HS:06}. It is known that $A$ is $K$-positive if and only if $A^{\top}$ is $K^*$-positive~\cite[Theorem 2.24]{AB-RJP:94}.
For a polyhedral cone $K$, the following lemma establishes a {\color{black}connection} between the $K$-Metzler and $K$-positive operators\footnote{ We note that, for a pointed and proper cone $K$, a proof for this Lemma can be found in~\cite[Theorem 8]{HS-MV:70}. Unfortunately, the proof in~\cite{HS-MV:70} does not generalize to the cones that are not proper or pointed.}.  We postpone the proof of this lemma to Appendix~\ref{app:B}. 
\begin{lemma}\label{thm:pos-mon}
             Let $A:\real^n\to \real^n$ be a linear operator and $K\subset \real^n$ be a polyhedral cone. The following statements are equivalent:
             \begin{enumerate}
             \item\label{p2:monotone} $A$ is $K$-Metzler,
             \item\label{p1:positive} there exists $h^*>0$ such that $I_n+hA$
               is $K$-positive for every $h\in [0,h^*]$,
               \item\label{p3:positive-diag} there exists $\alpha^*>0$ such that $A+\alpha I_n$
               is $K$-positive for every $\alpha\ge \alpha^*$.
               \end{enumerate}
             \end{lemma}

   \section{Gauge and dual gauge norm}
   
    In this section, we consider a proper cone $K\subseteq \real^n$ and we introduce the notion of gauge and dual gauge norms to define two different metric structures on $\real^n$. Moreover, we introduce the gauge matrix measure associated to a gauge norm. As we will see later, the gauge matrix measure plays an important role in contraction theory of $K$-monotone systems. Given a vector $\mathbf{e}\in \mathrm{int}(K)$, the gauge function (also called the Minkowski functional) $\|\cdot\|_{\mathbf{e},K}:\real^n\to \real_{\ge 0}$
 of the interval $[-\mathbf{e},\mathbf{e}]_K$ is defined as
 follows~\cite{CDA-RT:07}:
\begin{align}\label{eq:norm}
  \|v\|_{\mathbf{e},K} = \inf \setdef{\lambda\in \real_{\ge 0}}{v \in
  \lambda [-\mathbf{e},\mathbf{e}]_K}.
\end{align}
{\color{black}Using the identification $(\real^n)^* \cong \real^n$, we can define another metric structure on $\real^n$. If $K$ is pointed, for a vector $\mathbf{e}^*\in \mathrm{int}(K^*)$, the dual gauge function $\|\cdot\|^{\mathrm{d}}_{\mathbf{e}^*,K^*}:\real^n\to \real_{\ge 0}$
 of the interval $[-\mathbf{e}^*,\mathbf{e}^*]_{K^*}$ is defined as
 follows}:
\begin{align}\label{eq:dualnorm}
  \|v\|^{\mathrm{d}}_{\mathbf{e}^*,K^*} = \max \setdef{|\langle \eta, v\rangle|}{\eta\in [-\mathbf{e}^*,\mathbf{e}^*]_{K^*}}.
\end{align}
It is known that the gauge function defined in~\eqref{eq:norm} is a seminorm. For pointed proper cones, both the gauge function and the dual gauge function defined in~\eqref{eq:dualnorm} are norms.
\begin{proposition}[Gauge and dual gauge norms]\label{thm:norm}
Let $\real^n$ be a finite dimensional vector space and $K\subseteq \real^n$ be a proper cone. For every $\mathbf{e}\in \mathrm{int}(K)$, the following statements hold:
  \begin{enumerate}
      \item\label{p0:seminorm} the gauge function $\|\cdot\|_{\mathbf{e},K}$ is a seminorm on $\real^n$;
      \end{enumerate}
      {\color{black}Additionally, if $K\subseteq \real^n$ is pointed and $\mathbf{e}^*\in \mathrm{int}(K^*)$,} then
      \begin{enumerate}\setcounter{enumi}{1}
         \item\label{p1:gauge} the gauge function $\|\cdot\|_{\mathbf{e},K}$ is a norm on $\real^n$;
         \item\label{p2:dualgauge} the dual gauge function  $\|\cdot\|^{\mathrm{d}}_{\mathbf{e}^*,K^*}$ is a norm on $\real^n$.
      \end{enumerate}
      
\end{proposition}
  \smallskip
\begin{proof} 
Regarding part~\ref{p0:seminorm}, we first show that $\mathrm{int}([-\mathbf{e},\mathbf{e}]_K)$ is a neighborhood of the origin $0$.  Since $\mathbf{e}\in \mathrm{int}(K)$ and addition is a continuous function on $\real^n$, we get that $0\in \mathrm{int}(S)$ where $S = \setdef{x\in \real^n}{-\mathbf{e}\preceq_K x}$. Similarly, we have $-\mathbf{e}\in \mathrm{int}(-K)$. This implies that $0\in \mathrm{int}(S')$ where $S' = \setdef{x\in \real^n}{x\preceq_K \mathbf{e}}$. As a result, we get $0\in \mathrm{int}(S)\cap \mathrm{int}(S') = \mathrm{int}(S\cap S') = \mathrm{int}([-\mathbf{e},\mathbf{e}]_K)$. This means that $\mathrm{int}([-\mathbf{e},\mathbf{e}]_K)$ is a neighborhood of the origin. Moreover, every neighborhood of the origin is absorbent in the vector space $\real^n$~\cite[Theorem 4.3.6(b)]{LN-EB:85}. Therefore, by~\cite[Theorem 5.3.1]{LN-EB:85} the gauge function $\|\cdot\|_{\mathbf{e},K}$ is a seminorm on $\real^n$. Regarding part~\ref{p1:gauge}, since $K$ is pointed, it does not contain a non-trivial vector subspace. Now we can use~\cite[Exercise 5.105(d)]{LN-EB:85}, to show that the gauge functional $\|\cdot\|_{\mathbf{e},K}$ is a norm on $\real^n$. {\color{black}Regarding part~\ref{p2:dualgauge}, since $K$ is proper and pointed, the dual cone $K^*$ is proper and pointed too.} Now by part~\ref{p1:gauge}, one can define the norm $\|\cdot\|_{\mathbf{e}^*,K^*}$ on $\real^n$ by
\begin{align*}
 \|\phi\|_{\mathbf{e}^*,K^*} = \inf\setdef{\lambda\ge 0}{\phi\in \lambda [-\mathbf{e}^*,\mathbf{e}^*]_{K^*}}.   
\end{align*}
Then we have $\|v\|^{\mathrm{d}}_{\mathbf{e}^*,K^*} = \max \setdef{|\langle \eta, v\rangle|}{\eta\in [-\mathbf{e}^*,\mathbf{e}^*]_{K^*}}$ and, as a result, we have $\|v\|^{\mathrm{d}}_{\mathbf{e}^*,K^*} = \max \setdef{|\langle \eta, v\rangle|}{\|\eta\|_{\mathbf{e}^*,K^*}\le 1}$. Thus, $\|\cdot\|^{\mathrm{d}}_{\mathbf{e}^*,K^*}$ is the dual norm to $\|\cdot\|_{\mathbf{e}^*,K^*}$ on $\real^n$.\end{proof}
\smallskip

For a polyhedral cone $K\subseteq \real^n$ with a representation $(H,V)$, there exist closed-form expressions for the gauge and the dual gauge norm.
\begin{lemma}[Formula for the gauge seminorms]\label{lem:gaugenorm}
Suppose that $K\subset\real^n$ is a proper polyhedral cone with a representation $(H,V)$ such that $H\in \real^{m\times n}$ and $V\in \real^{n\times q}$ and with $\mathbf{e}\in \mathrm{int}(K)$. Then 
\begin{enumerate}
    \item\label{p1:gaugenorm} $\|x\|_{\mathbf{e},K}=\|\mathrm{diag}(H\mathbf{e})^{-1}Hx\|_{\infty}$.
    \end{enumerate}
    {\color{black}Additionally, if $K\subset \real^n$ is pointed and $\mathbf{e}^*\in \mathrm{int}(K^*)$,} then
    \begin{enumerate}\setcounter{enumi}{1}
    \item\label{p2:dualguagenorm} $\|x\|^{\mathrm{d}}_{\mathbf{e}^*,K^*}=\|\mathrm{diag}(V^{\top}\mathbf{e}^*)V^{\dagger}x\|_1$.
\end{enumerate}
\end{lemma}
  \smallskip
\begin{proof}
Regarding part~\ref{p1:gaugenorm}, note that, by definition of the gauge norm,  we have $\|v\|_{\mathbf{e},K} = \inf\setdef{\lambda}{ -\lambda\mathbf{e}\preceq_K v \preceq_K \lambda\mathbf{e}}$. Since $\mathbf{e}\in \mathrm{int}(K)$, we get that $H\mathbf{e}>\vect{0}_m$~\cite[Proposition 1.1]{SW:01}. Using Lemma~\ref{thm:inequality}, we get $\|v\|_{\mathbf{e},K} = \inf\setdef{\lambda}{-\lambda H\mathbf{e}\le Hv \le \lambda H\mathbf{e}}$. 
Multiplying the above inequalities $-\lambda H\mathbf{e}\le Hv \le \lambda H\mathbf{e}$ by the positive diagonal matrix $\mathrm{diag}(H\mathbf{e})^{-1}$, we get 
\begin{align*}
\|&v\|_{\mathbf{e},K} = \inf\setdef{\lambda}{-\lambda \vect{1}_m\le \mathrm{diag}(H\mathbf{e})^{-1} Hv \le \lambda \vect{1}_m} \\ & = \inf\setdef{\lambda}{\|\mathrm{diag}(H\mathbf{e})^{-1}Hv\|_{\infty}\le \lambda} = \|\mathrm{diag}(H\mathbf{e})^{-1}Hv\|_{\infty}.
\end{align*}
Regarding part~\ref{p2:dualguagenorm}, note that 
\begin{align*}
    \|v\|^{\mathrm{d}}_{\mathbf{e}^*,K^*} = \max\setdef{|\langle\xi, v\rangle|}{-\mathbf{e}^* \preceq_{K^*} \xi \preceq_{K^*} \mathbf{e}^*}.
\end{align*}
 Using Lemma~\ref{thm:inequality}, we get 
 \begin{align*}
    \|v\|^{\mathrm{d}}_{\mathbf{e}^*,K^*} = \max\setdef{|\langle\xi, v\rangle|}{-V^{\top}\mathbf{e}^* \le V^{\top}\xi \le V^{\top}\mathbf{e}^*}.
\end{align*}
Since $\mathbf{e}^*\in \mathrm{int}(K^*)$, we have $V^{\top}\mathbf{e}^*>\vect{0}_q$~\cite[Proposition 1.1]{SW:01}. Multiplying the above inequalities by the positive diagonal matrix $\mathrm{diag}(V^{\top}\mathbf{e}^*)^{-1}$, we get
 \begin{align*}
    \|v\|^{\mathrm{d}}_{\mathbf{e}^*,K^*} = \max\setdef{|\langle\xi, v\rangle|}{-\vect{1}_q \le \mathrm{diag}(V^{\top}\mathbf{e}^*)^{-1} V^{\top}\xi \le \vect{1}_q}.
\end{align*}
Since $K$ is proper, we have $\mathrm{int}(K)\ne \emptyset$ and $V$ has full row rank. As a result, we get  $(V^{\top})^{\dagger}V^{\top} = I_n$ and $(V^{\top})^{\dagger} = (V^{\dagger})^{\top}$. Thus,
\begin{align*}
|\langle\xi, v\rangle| & = |\langle(V^{\dagger})^{\top}V^{\top} \xi, v\rangle| = |\langle V^{\top} \xi, V^{\dagger}v\rangle| \\ & = |\langle \mathrm{diag}(V^{\top}\mathbf{e}^*)^{-1}V^{\top} \xi, \mathrm{diag}(V^{\top}\mathbf{e}^*) V^{\dagger}v\rangle|.    
\end{align*}
As a result, $\|v\|^{\mathrm{d}}_{\mathbf{e}^*,K^*} = \max\setdef{|\langle \eta, \mathrm{diag}(V^{\top}\mathbf{e}^*)V^{\dagger}v\rangle|}{\|\eta\|_{\infty}\le 1}  = \|\mathrm{diag}(V^{\top}\mathbf{e}^*)V^{\dagger}v\|_{1}$, 
where the last equality holds because the $\ell_1$-norm is the dual of the $\ell_{\infty}$-norm on $\real^q$.
\end{proof}
\smallskip
Given a proper polyhedral cone $K\subseteq \real^n$ with $\mathbf{e}\in \mathrm{int}(K)$, the matrix semi-measure associated to the gauge seminorm $\|\cdot\|_{\mathbf{e},K}$ is denoted by $\mu_{\mathbf{e},K}$. Note that $\mu_{\mathbf{e},K}$ is a matrix measure if and only if $K$ is pointed. Now, we present two examples of polyhedral cones and their associated gauge and dual gauge norms. 
    \paragraph*{Standard Euclidean cone}
    The set of all non-negative vectors $\real^n_{\ge 0}$ is a pointed proper cone in $\real^n$
    with a non-empty interior. The partial order
     associated with $\real^n_{\ge 0}$ is the standard component-wise order on
     $\real^n$, i.e., $x\le y$ if we have $x_i\le y_i$, for every
     $i\in\{1,\ldots,n\}$. For $\mathbf{e}=\vect{1}_n$, the gauge norm
     $\|\cdot\|_{\vect{1}_n,\real^n_{\ge 0}}$ is the standard
     $\ell_{\infty}$-norm on $\real^n$. It can be shown that $K^*=\real^{n}_{\ge 0}$ and, by choosing $\mathbf{e}^*=\vect{1}_n$, the dual gauge norm $\|\cdot\|^{\mathrm{d}}_{\vect{1}_n,\real^{n}_{\ge 0}}$ is the standard $\ell_1$-norm on $\real^n$;
       \smallskip
     \paragraph*{$1$-norm cone} For every $S\subseteq \{2,3,\ldots,n\}$, we define the linear functional $\phi_S:\real^n\to \real$ as follows:
  \begin{align*}
    \phi_S(v) = v_1 + \sum_{j\in S}v_j - \sum_{k\not\in S\cup\{1\}}v_k
  \end{align*}
  and we define $K\subseteq \real^n$ as the pointed proper polyhedral cone generated by
  $\{\phi_S\}$ for every $S \subseteq \{2,3,\ldots,n\}$, i.e.,
\begin{align*}
  K =\setdef{x\in \real^n}{\langle\phi_S,x\rangle\ge 0, \;\;\; S\subseteq \{2,3,\ldots,n\}}.
\end{align*}
By choosing $\mathbf{e}=(1\;  0\; \cdots \;  0)^{\top}\in \real^n$, we get $\|v\|_{\mathbf{e},K} = \|v\|_1$. 

The next theorem provides a characterization of the gauge and dual gauge matrix measures for $K$-Metzler operators. This result can be considered as a generalization of the classical formulas for $\ell_1$- and the $\ell_{\infty}$-matrix measures for Metzler matrices~\cite[Equations (4) and (5)]{SC:19}. We postpone the proof of this theorem to Appendix~\ref{ap:C}. 
\smallskip
             \begin{theorem}[Characterization of the gauge matrix measures]\label{thm:contraction}
               Consider a proper polyhedral cone $K\subseteq \real^n$ with a representation $(H,V)$ and with $\mathbf{e}\in \mathrm{int}(K)$. Suppose that $A:\real^n\to \real^n$ is a $K$-Metzler linear operator. The following statements are equivalent: 
               \begin{enumerate}
               \item\label{p1} $\mu_{\mathbf{e},K}(A)\le c$,
            \smallskip
               \item\label{p2} $A\mathbf{e}\preceq_K c\mathbf{e}$, 
               \smallskip
               \item\label{p3} $HA\mathbf{e}\le c H\mathbf{e}$. 
               \end{enumerate}
               {\color{black}Additionally, if $K$ is pointed and $\mathbf{e}^*\in \mathrm{int}(K^*)$,} the following statements are equivalent:
               \begin{enumerate}\setcounter{enumi}{3}
               \item\label{p1:dual} $\mu^{\mathrm{d}}_{\mathbf{e}^*,K^*}(A)\le c$,
                \smallskip
               \item\label{p2:dual} $A^{\top}\mathbf{e}^*\preceq_{K^*} c\mathbf{e}^*$,
               \smallskip
               \item\label{p3:dual} $V^{\top} A^{\top}\mathbf{e}^* \le c V^{\top}\mathbf{e}^*$.
               \end{enumerate}
             \end{theorem}
\section{Monotone system on polyhedral cones}
In this section, we study monotonicity of a control system with respect to polyhedral cones. Consider the dynamical system 
\begin{align}\label{eq:dynamicalsystem}
    \dot{x} = f(x,u),
\end{align}
where $x\in \real^n$ is the system state and $u\in \real^p$ is the control input. We assume that $K\subset \real^n$ is a cone. 

  \begin{definition}[Monotone systems] Consider the control system~\eqref{eq:dynamicalsystem} with a cone $K\subseteq \real^n$. Then the system~\eqref{eq:dynamicalsystem} is $K$-monotone if, for every $x_0\preceq_K y_0$ and every $u\in \real^p$, we have
  \begin{align*}
      x_u(t) \preceq_K y_u(t),\quad \mbox{ for every }t\in \real_{\ge 0}, 
  \end{align*}
where $t\mapsto x_u(t)$ and $t\mapsto y_u(t)$ are trajectories of~\eqref{eq:dynamicalsystem} with {\color{black} input} $u$ starting from $x_0$ and $y_0$, respectively. 
     \end{definition}
\smallskip
One can show that the control system~\eqref{eq:dynamicalsystem} is $K$-monotone if and only if for every $x\preceq_K y$, every $\phi\in K^*$ satisfying $\langle
       \phi, x\rangle = \langle \phi, y\rangle$, and every $u\in \real^p$, we have~\cite[Theorem 1]{DA-EDS:03}
     \begin{align*}
       \langle \phi, f(x,u)\rangle \le  \langle \phi, f(y,u)\rangle.
       \end{align*}
When the map $x\mapsto f(x,u)$ is continuously differentiable and $K$ is a proper cone, one can show that the control system~\eqref{eq:dynamicalsystem} is $K$-monotone if and only if $D_xf(x,u)$ is $K$-Metzler, for every $(x,u)\in \real^n\times \real^p$~\cite[Theorem 3.5]{MWH-HS:06}.
Our first result of this section provides three equivalent characterizations for the control system~\eqref{eq:dynamicalsystem} to be $K$-monotone with respect to a proper polyhedral cone $K$ with a representation $(H,V)$. 
             \begin{theorem}[Characterization of monotonicity]\label{thm:monotone}
             Consider the control system~\eqref{eq:dynamicalsystem} with continuously differentiable $f$. Let $K\subseteq \real^n$ be a proper polyhedral cone with a representation $(H,V)$ such that $H\in \real^{m\times n}$ and $V\in \real^{n\times q}$. The following statements are equivalent:
               \begin{enumerate}
               \item\label{p1:mon} the dynamical system~\eqref{eq:dynamicalsystem} is $K$-monotone;
               \item\label{p2:mon} {\color{black}there exists $\alpha^*:\real^n\times \real^p\to \real$ such that
               \begin{align}\label{eq:cond1}
                 H(D_xf(x,u) + \alpha I_n) V\ge \vect{0}_{m\times q},
               \end{align}
               for every $(x,u)\in \real^n\times \real^p$ and every $\alpha\ge \alpha^*(x,u)$}; 
               \item\label{p3:mon} there exists $P:\real^{n}\times \real^p\to \mathbb{M}_m$ such that
               \begin{align}\label{eq:cond2}
                HD_xf(x,u) = P(x,u) H,
               \end{align}
               for every $(x,u)\in \real^n\times \real^p$; 
               \item\label{p4:mon} there exists $Q:\real^{n}\times \real^{p}\to \mathbb{M}_q$ such that
               \begin{align}\label{eq:cond3}
                   D_xf(x,u)V = VQ(x,u),
               \end{align}
               for every $(x,u)\in \real^n\times \real^p$.
               \end{enumerate}
             \end{theorem}
               \smallskip
             \begin{proof}
             First note that, using Lemma~\ref{thm:pos-mon}, the control system~\eqref{eq:dynamicalsystem} is $K$-monotone if and only if there exists $\alpha^*:\real^n\times \real^p\to \real_{\ge 0}$ such that $D_xf(x,u)+\alpha I_n$ is $K$-positive for every $(x,u)\in \real^n\times \real^p$ and every $\alpha\ge \alpha^*(x,u)$. Regarding \ref{p1:mon} $\iff$ \ref{p2:mon}, the result
               then follows from~\cite[Theorem 4.1]{FB-MF-EH:74}. Regarding \ref{p1:mon} $\iff$ \ref{p3:mon}, we denote the $i$th row of the matrix $H$ by $h_i \in \real^{1\times n}$, for every $i\in \{1,\ldots,m\}$. By the equivalence of~\ref{p1:mon} and~\ref{p2:mon}, the control system~\eqref{eq:dynamicalsystem} is $K$-monotone if and only if there exists $\alpha^*:\real^n\times \real^p\to \real_{\ge 0}$ such that, for every $(x,u)\in \real^n\times \real^p$, every $\alpha>\alpha^*(x,u)$, and every $i\in\{1,\ldots,m\}$, we have $h_i(D_xf(x,u) + \alpha I_n)v\ge 0$ for every $v\in \real^n$ such that $Hv\ge \vect{0}_m$. Now, using Farkas's Lemma~\cite[Proposition 1.8]{GMZ:12}, the control system~\eqref{eq:dynamicalsystem} is $K$-monotone if and only if there exists $\eta_i\ge \vect{0}_m$ such that 
               $h_i(D_xf(x,u) + \alpha I_n)= \eta^{\top}_i H$. Therefore, the control system~\eqref{eq:dynamicalsystem} is $K$-monotone if and only if, for every $(x,u)\in \real^n\times \real^p$, we have $H D_xf(x,u) = (Q(x,u) - \alpha I_m )H$, for some non-negative matrix $Q(x,u)\ge \vect{0}_{m\times m}$. The result follows by defining $P(x,u) := Q(x,u) - \alpha I_m$ and noting that $P(x,u)$ is a Metzler matrix, for every $(x,u)\in \real^n\times \real^p$.

               Regarding \ref{p1:mon} $\iff$ \ref{p4:mon}, using the equivalence of~\ref{p1:mon} and~\ref{p2:mon}, one can easily show that the control system~\eqref{eq:dynamicalsystem} is $K$-monotone if and only if, for every $(x,u)\in \real^n\times \real^p$ and every $\eta\in\real^q_{\ge 0}$, there exists $\xi\in \real^q_{\ge 0}$ such that $(D_xf(x,u)+\alpha(x,u) I_n)V\eta = V\xi$. In turn, the last statement is equivalent to the following sentence: for every $(x,u)\in \real^n\times \real^p$, $(D_xf(x,u)+\alpha(x,u) I_n)V = VP(x,u)$, for some non-negative matrix $P(x,u)\ge \vect{0}_{q\times q}$. The result then follows by defining the Metzler matrix $Q(x,u)\in \mathbb{M}_q$ by $Q(x,u)=P(x,u)-\alpha(x,u)I_q$, for every $(x,u)\in \real^n\times \real^p$. 
               \end{proof}
                 \smallskip
             \begin{remark}
             The following remarks are in order. 
             \begin{enumerate}
                 \item In~\cite[Theorem 4.3]{YK-FF:22}, the equivalence of $K$-monotonicty of the system~\eqref{eq:dynamicalsystem} and condition~\eqref{eq:cond2} is shown for proper pointed cones. By appealing to the novel Lemma~\ref{thm:pos-mon}, Theorem~\ref{thm:monotone} extends this result to general polyhedral cones $K$ without assuming that $K$ is proper or pointed.
                 {\color{black}\item In~\cite[Lemma 4]{YO-HM-SK:84}, the equivalence of $K$-monotonicity and condition~\eqref{eq:cond3} is shown for linear systems and for arbitrary polyhedral cone $K$. Theorem~\ref{thm:monotone} generalizes this result to nonlinear systems.}
                 \item One can compare the conditions~\eqref{eq:cond1},~\eqref{eq:cond2}, and~\eqref{eq:cond3} in terms of their computational complexity. Checking $K$-monotonicity using Theorem~\ref{thm:monotone}\ref{p2:mon} requires knowledge of both $H$-rep and the $V$-rep of the polyhedral cone $K$ but one only needs to solve condition~\eqref{eq:cond1} for a scalar $\alpha^*(x,u)\in \real$. On the other hand, checking $K$-monotonicity using the Theorem~\ref{thm:monotone}\ref{p3:mon}  (resp. Theorem~\ref{thm:monotone}\ref{p4:mon}) only requires the knowledge of the $H$-rep (resp. $V$-rep) of the polyhedral cone $K$ but one needs to solve condition~\eqref{eq:cond2}  (resp. condition~\eqref{eq:cond3}) for an $m\times m$ Metzler matrix $P(x,u)\in \mathbb{M}_m$ (resp. a $q\times q$ Metzler matrix $Q(x,u)\in\mathbb{M}_q$).
             \end{enumerate}
               \end{remark}
               \smallskip

            Now, we can state our main result which characterizes contractivity of a $K$-monotone system with respect to the gauge norm. Note that, if the polyhedral cone $K$ is pointed, then a similar result can be obtained for contractivity of the $K$-monotone system  with respect to the dual gauge norm. We omit this result for brevity of presentation. 
            
              \smallskip
               
               \begin{theorem}[Semi-contraction for the gauge seminorm]\label{thm:contraction-dualnom}
               Consider the control system~\eqref{eq:dynamicalsystem}. Let $K\subseteq \real^n$ be a proper polyhedral cone with a  representation $(H,V)$. Let $\mathbf{e}\in \mathrm{int}(K)$, $c\in \real$, and $\|\cdot\|_{\mathcal{U}}$ be a norm on $\real^p$. Suppose that the control system~\eqref{eq:dynamicalsystem} is $K$-monotone and $D_x f(x,u)\mathrm{Ker}(H)\subseteq \mathrm{Ker}(H)$, for every $(x,u)\in \real^n\times \real^p$. The following statements are equivalent:
               \begin{enumerate}
                   \item\label{p1:contraction} $H D_xf(x,u) \mathbf{e} \le -c H\mathbf{e}$, for every $(x,u)\in \real^n\times \real^p$;
                     \smallskip
                   \item\label{p2:contraction} any two trajectories $x_u(t),y_u(t)$ with the same continuous input signal $u:\real_{\ge 0}\to \real^p$ satisfy:
                   \begin{align*}
                       \|x_u(t)-y_u(t)\|_{\mathbf{e},K}\le e^{-ct}\|x_u(0)-y_u(0)\|_{\mathbf{e},K}. 
                   \end{align*}
                   \item\label{p3:robustness} any two trajectories $x_u(t),y_v(t)$ with different continuous input signals $u,v:\real_{\ge 0}\to \real^p$ satisfy:
                   \begin{multline*}
                       \|x_u(t)-y_v(t)\|_{\mathbf{e},K}\le e^{-ct}\|x_u(0)-y_v(0)\|_{\mathbf{e},K} \\ + \frac{\ell(1-e^{-ct})}{c}\sup_{\tau\in [0,t]}\|u(\tau)-v(\tau)\|_{\mathcal{U}},
                   \end{multline*}
               \end{enumerate}
               where $\ell=\sup_{x,u\in \real^n\times \real^p} \sup_{\eta\in \real^p} \frac{\|D_uf(x,u)\eta\|_{\mathbf{e},K}}{\|\eta\|_{\mathcal{U}}}$. Moreover, if the cone $K$ is pointed, condition~\ref{p1:contraction} holds for some $c>0$, and $u\in \real^p$ is a constant input signal,  then \begin{enumerate}\setcounter{enumi}{3}
                  \item\label{p3:eqpt} the system~\eqref{eq:dynamicalsystem} has a unique globally exponentially stable equilibrium point $x^*\in \real^n$;
                  \item\label{p4:lyapunov} the functions 
                  \begin{align*}
                       V_1(x) =& \|\mathrm{diag}(H\mathbf{e})^{-1}H(x-x^*)\|_{\infty}, \\ V_2(x) =& \|\mathrm{diag}(H\mathbf{e})^{-1}Hf(x,u)\|_{\infty}.
                   \end{align*}
                   are global Lyapunov functions for~\eqref{eq:dynamicalsystem}. 
              \end{enumerate}
               \end{theorem}
                 \smallskip
               \begin{proof}
                Regarding~\ref{p1:contraction} $\iff$~\ref{p2:contraction}, first note that the control system~\eqref{eq:dynamicalsystem} is $K$-monotone and thus $D_xf(x,u)$ is $K$-Metzler for every $(x,u)\in \real^{n}\times \real^p$~\cite[Theorem 3.5]{MWH-HS:06}. Now, using Theorem~\ref{thm:contraction}, the condition in part~\ref{p1:contraction} is equivalent to $\mu_{\mathbf{e},K}(D_x f(x,u))\le -c$, for every $(x,u)\in \real^n\times \real^p$. The result then follows from~\cite[Theorem 59]{AD-SJ-FB:20o}. Regarding~\ref{p1:contraction} $\implies$~\ref{p3:robustness}, by Lemma~\ref{lem:gaugenorm}\ref{p1:gaugenorm}, for every $x\in \real^n$, we have $\|x\|_{\mathbf{e},K} = \|\mathrm{diag}(H\mathbf{e})^{-1}Hx\|_{\infty}$. The result is then straightforward by replacing the norms in the proof of~\cite[Theorem 37(ii)]{AD-SJ-FB:20o} by the seminorm $x\mapsto \|\mathrm{diag}(H\mathbf{e})^{-1}Hx\|_{\infty}$.
                Regarding~\ref{p3:robustness} $\implies$~\ref{p2:contraction}, the result is straightforward by choosing $v(t)=u(t)$, for every $t\in \real_{\ge 0}$.
                Regarding \ref{p1:contraction} $\implies$~\ref{p3:eqpt} and \ref{p4:lyapunov}, we note that if $K$ a is proper pointed cone, by Proposition~\ref{thm:norm}, the gauge function $\|\cdot\|_{\mathbf{e},K}$ is a norm. The results then follow by~\cite[Theorem 3.8]{FB:22-CTDS}.
                \end{proof}
                  \smallskip
               \begin{remark}[Comparison with the literature]\;\;
               \begin{enumerate}
                   \item Theorem~\ref{thm:contraction-dualnom} is a generalization of \cite[Theorem 2]{SC:19} to $K$-monotone systems for a proper polyhedral cone $K$. Additionally, Theorem~\ref{thm:contraction-dualnom} provides an incremental input-to-state robustness bound for contractive $K$-monotone systems.    
                   \item In~\cite[Theorem 4.5]{YK-FF:22}, a sufficient condition for exponential incremental  stability of a $K$-monotone system is proposed based upon embedding the system into a higher dimensional space. In comparison, our Theorem~\ref{thm:contraction-dualnom} presents a necessary and sufficient condition for contractivity of $K$-monotone systems with respect to the gauge norm $\|\cdot\|_{\mathbf{e},K}$. It is worth mentioning that exponential incremental stability is a weaker condition than contractivity with respect to any norm. 
                   \item Given a proper polyhedral cone $K\subseteq \real^n$, the sufficient condition for exponential incremental stability in~\cite[Theorem 4.5]{YK-FF:22} requires searching for a vector $v\in\real^m$ and a Metzler matrix $P\in \real^{m\times m}$. However, the condition in Theorem~\ref{thm:contraction-dualnom}\ref{p1:contraction} only requires searching for one scalar, i.e., {\color{black}$c\in \real$}. Thus, in cases when the proper polyhedral cone $K$ is given, the sufficient condition in Theorem~\ref{thm:contraction-dualnom}\ref{p1:contraction} is computationally more efficient than the condition in~\cite[Theorem 4.5]{YK-FF:22}. 
                   
                   \item For a control system on $\real^n$ with a globally stable equilibrium point $x^*$, the search for a quadratic Lypaunov function of the form $V(x) = (x-x^*)^{\top}\Sigma (x-x^*)$ requires solving for $\frac{n(n-1)}{2}$ entries of the positive definite matrix $\Sigma\in \real^{n\times n}$. In this context, Theorem~\ref{thm:contraction-dualnom}\ref{p1:contraction} provides a scalable approach to construct global polygonal Lyapunov functions for the $K$-monotone system by searching for $n$ components of a vector $\mathbf{e}\in \mathrm{int}(K)$.

               \end{enumerate}
               \end{remark}

\section{Applications}

In this section, we present two applications of our framework for analysis and design of systems. 


\subsection{Monotone edge flows in dynamic networks}
Consider a network of interconnected compartments, where the state of the compartment $i$ is described by $x_i\in \real$, for every $i\in \{1,\ldots,n\}$. The interconnection of the compartments is described by a connected undirected graph $G=(\mathcal{V},\mathcal{E})$, where $\mathcal{V}=\{1,\ldots,n\}$ is the set of nodes and $\mathcal{E}\subseteq \mathcal{V}\times \mathcal{V}$ is the set of edges described by unordered pairs of vertices. The dynamics of the network is given by 
  \begin{align}\label{eq:flowvector}
    \dot{x}=f(x)
\end{align}
where $x = (x_1,\ldots,x_n)^{\top}$ and  $f(x) = (f_1(x),\ldots, f_n(x))^{\top}$.  We assume that the vector field $f$ satisfies the following translation-invariance property:
\begin{align}\label{eq:translationinv}
    f(x+ c\vect{1}_n) = f(x), \quad\mbox{ for every }x\in \real^n, \;\; c\in \real.
\end{align}
In the literature, the translation-invariance property~\eqref{eq:translationinv} is sometimes referred to as \emph{plus-homogeneity}~\cite{DD-Mf-AG:23} and has shown to play a critical role in stability of $K$-monotone dynamical systems~\cite{DA-ES:08}.

For every edge $e=(i,j)\in \mathcal{E}$, the edge flow from compartment $i$ to compartment $j$ is defined by $x_i-x_j$. Given an edge orientation for the graph $G$, the vector of the flows is given by $B^{\top}x\in \real^m$, where $B\in \real^{n\times m}$ is the incidence matrix of the graph $G$ associated to the given edge orientation~\cite[Chapter 9]{FB:22}. For many real-world interconnected systems, including power grids and traffic networks, edge flows correspond to physical quantities and play a crucial role in safety and security analysis of networks. In this section, we study the evolution of the edge flows in the dynamic flow network~\eqref{eq:flowvector}. We start our analysis of the flow dynamics~\eqref{eq:flowvector} by defining the set
\begin{align*}
    \mathcal{I}_{\mathcal{E}}:=\setdef{v\in \real^n}{v_i\ne v_j, \mbox{ for every }(i,j)\in \mathcal{E}}.
\end{align*}
Given $v\in \mathcal{I}_{\mathcal{E}}$, one can assign an edge orientation to the graph $G=(\mathcal{V},\mathcal{E})$ such that if $B_v\in \real^{n\times m}$ is the incidence matrix of $G$ associated to this edge orientation, we have $B^{\top}_v v> \vect{0}_m$. For every $v\in \mathcal{I}_{\mathcal{E}}$, we define the cone $K^v_G\subset\real^n$ by
\begin{align}\label{eq:incidencecone}
    K^v_G =\setdef{x\in \real^n}{B_v^{\top}x\ge \vect{0}_m}.
\end{align}
By definition of the incidence matrix $B_v$, we have $B^{\top}_v v > \vect{0}_m$. {\color{black}This implies that $v\in \mathrm{int}(K^v_G)$}. Therefore $K^v_{G}$ is a proper cone. However, $K^v_G$ is not a pointed cone. This is because  $B_v^{\top}\vect{1}_n=\vect{0}_m$ and thus $\mathrm{span}\{\vect{1}_n\}\in K^v_G\cap (-K^v_G)$. The next theorem studies monotonicity and contractivity of edge flows in dynamic flow networks.

\begin{theorem}[Monotonicity and contractivity of edge flows]\label{thm:edgeflow}
Consider the dynamic flow network~\eqref{eq:flowvector} over an undirected connected graph $G=(\mathcal{V},\mathcal{E})$ and let $v\in \mathcal{I}_{\mathcal{E}}$.  
The following statements are equivalent:
\begin{enumerate}
    \item\label{p1:flow} for every $x\in \real^n$, there exists $P(x)\in \mathbb{M}_m$ such that
    \begin{align*}
        B_v^{\top} D_xf(x) = P(x) B_v^{\top};
    \end{align*}
    \item\label{p2:flow} for every two trajectories $x(t),y(t)$ of the system~\eqref{eq:flowvector} satisfying $B_v^{\top}x(0)\le B_v^{\top}y(0)$, we have
\begin{align*}
    B_v^{\top}x(t) \le B_v^{\top}y(t), \quad\mbox{ for every }t\in \real_{\ge 0}.
\end{align*}
\end{enumerate}
Additionally, if condition~\ref{p1:flow} holds, $c\in \real$, and $\mathbf{e}\in \mathrm{int}(K^v_G)$, then the following statements are equivalent:
\begin{enumerate}\setcounter{enumi}{2}
\item \label{p3:flow} $B_v^{\top}D_xf(x)\mathbf{e}\le -c B_v^{\top}\mathbf{e}$, for every $x\in \real^n$.
\smallskip
    \item\label{p4:flow} for every two trajectories $x(t),y(t)$ of the system~\eqref{eq:flowvector},
\begin{multline*}
   \big\|\mathrm{diag}(B_v^{\top}\mathbf{e})^{-1}\big(B_v^{\top}x(t)-B_v^\top y(t)\big)\big\|_{\infty} \\ \le e^{-ct} \big\|\mathrm{diag}(B_v^{\top}\mathbf{e})^{-1}\big(B_v^{\top}x(0)-B_v^\top y(0)\big)\big\|_{\infty},
\end{multline*}
for every $t\in \real_{\ge 0}$.
\end{enumerate}
\end{theorem}
\smallskip
\begin{proof}
The equivalence \ref{p1:flow}~$\iff$~\ref{p2:flow} follows from Theorem~\ref{thm:monotone}\ref{p2:mon} and Lemma~\ref{thm:inequality} applied to the cone $K^v_G$. Regarding the equivalence \ref{p3:flow}~$\iff$~\ref{p4:flow}, first note that by translation-invariance law~\eqref{eq:translationinv}, we have $ D_xf(x) \vect{1}_n = \vect{0}_n$. On the other hand, the graph $G$ is connected and thus $\mathrm{Ker}(B^{\top}_v)=\mathrm{span}(\vect{1}_n)$. As a result, we have $D_xf(x)\mathrm{Ker}(B_v^{\top}) = \{\vect{0}_n\}\subset \mathrm{Ker}(B_v^{\top})$. The result then follows from Theorem~\ref{thm:contraction-dualnom} and Lemma~\ref{lem:gaugenorm}\ref{p1:gaugenorm}.
\end{proof}

\begin{example}[Edge flows in averaging systems]\label{ex:ave}
Let $G =(\mathcal{V},\mathcal{E})$ be a network with an undirected connected graph shown in Figure~\ref{fig:cycle} and consider the following continuous-time averaging system on $G$:
\begin{align}\label{eq:ave}
    \dot{x}=-Lx, 
\end{align}
where $L=\left[\begin{smallmatrix}2 & -1 & -1 & 0\\ -1 & 2 & -1 & 0\\ -1 & -1 & 3 & -1\\ 0 &0 &-1 &1\end{smallmatrix}\right]$ is the Laplacian matrix of $G$. 
\begin{figure}
 \begin{center}
		\includegraphics[width =0.6\linewidth,clip]{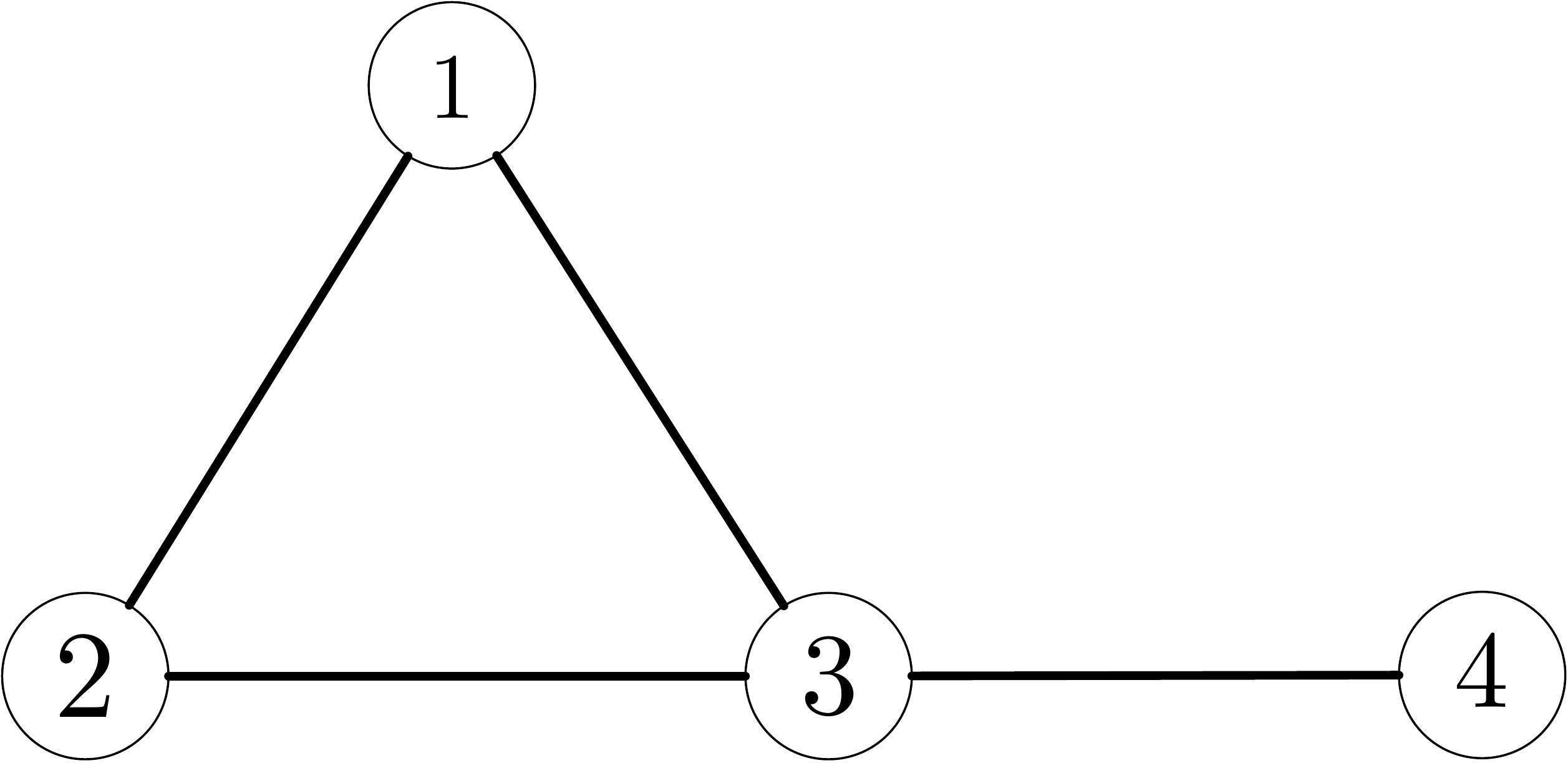}
	\end{center}
	\caption{The structure of the graph $G$ in Example~\ref{ex:ave}}
	\label{fig:cycle}
\end{figure}
 Let $v=(0,1,2,3)^{\top}\in \mathcal{I}_{\mathcal{E}}$. Then one can see that $B_v=\left[\begin{smallmatrix}1 & 0 & 1 & 0\\ -1 & 1 & 0 & 0\\ 0 & -1 & -1 & 1\\ 0 &0 & 0 &-1\end{smallmatrix}\right]$ and $B_v$ is the incidence matrix of $G$ associated with the orientation shown in Figure~\ref{fig:cycle}. We note that $L=B_vB_v^{\top}$. One can check that 
\begin{align*}
    -B_v^{\top}L = \left[\begin{smallmatrix} -3   &  3 &    0  &  0\\
     0 & -3 & 4 & -1\\
    -3 & 0 & 4 & -1\\
     1 & 1 & -4 & 2\end{smallmatrix}\right] = \left[\begin{smallmatrix} -3   &  0 &    0  &  0\\
     0 & -3 & 0 & 1\\
    0 & 0 & -3 &1\\
     0 & 1 & 1 & -2\end{smallmatrix}\right] B_v^{\top}.
\end{align*}
Since $\left[\begin{smallmatrix} -3   &  0 &    0  &  0\\
     0 & -3 & 0 & 1\\
    0 & 0 & -3 &1\\
     0 & 1 & 1 & -2\end{smallmatrix}\right]\in \mathbb{M}_4$, by Theorem~\ref{thm:edgeflow}, the edge flows of the averaging system~\eqref{eq:ave} are monotone. Moreover, one can pick $\mathbf{e}=(1.5,\;\; 1.4,\;\; 1,\;\; 0.1)^{\top}$ and check that $B_v^{\top}\mathbf{e} = (0.1,\;\; 0.4,\;\; 0.5,\;\; 0.9)^{\top}>\vect{0}_{4}$. Thus, $\mathbf{e}\in \mathrm{int}(K)$ and, we have 
     \begin{align*}
         -B_v^{\top}L\mathbf{e} = \left[\begin{smallmatrix} -0.3 \\-0.3\\-0.6\\-0.9 \end{smallmatrix}\right] \le - \frac{3}{4} B_v^{\top}\mathbf{e}. 
     \end{align*}
     Therefore, by Theorem~\eqref{thm:edgeflow}, the edge flows of the averaging systems~\eqref{eq:ave} are contracting with rate $c=\frac{3}{4}$. Alternatively, one can define the edge flow variable $z=B_v^{\top}x$ to get the edge flow dynamics:
\begin{align}\label{eq:edgeflowdynamic}
    \dot{z} = B_v^{\top}\dot{x} = - B_v^{\top}B_v z = -L_{\mathcal{E}}z = \left[\begin{smallmatrix}-2 & 1 & -1 & 0\\ 1 & -2 & -1 & 0\\ -1 & -1 & -2 & 1 \\0 & 1 & 1& -2 \end{smallmatrix}\right] z.
\end{align}
The matrix $L_{\mathcal{E}}=B_v^{\top}B_v$ is called the edge Laplacian of $G$. It is interesting to note that the edge Laplacian matrix $L_{\mathcal{E}}$ is not Metzler and $\lambda_{\max}(-L_{\mathcal{E}})=0$. Thus, one cannot deduce monotonicity or contractivity of the edge flows using the edge flow dynamics~\eqref{eq:edgeflowdynamic}. 
\end{example}

\subsection{Scalable control design with safety guarantees}

Monotone system theory has been successfully used for scalable control design in cooperative systems~\cite{AR:15} and in systems with rectangular safety constraints~\cite{EMV-PM:14}. However, in many applications, due to the nature of the problem, estimating the safe set using hyper-rectangles can either make the control design infeasible or can lead to overly-conservative results. In this subsection, we develop a scalable approach for state feedback design with safety guarantees, where we under-approximate the safe set using polytopes. Consider the following control system:
\begin{align}\label{eq:safety}
    \dot{x} = f(x) + Bu + Cw
\end{align}
where $x\in \real^n$ is the state of the system, $u\in \real^p$ is the control input, and $w\in \real^r$ is the vector of disturbance. We assume that $B\in \real^{n\times p}$ and $C\in \real^{n\times r}$ and the disturbance is bounded, i.e., there exists $\underline{w}\le \overline{w}$ such that $\underline{w}\le w(t)\le \overline{w}$, for every $t\in \real_{\ge 0}$. We assume that the origin $\vect{0}_n$ is a (possibly unstable) equilibrium point of $f$ and there exists a finite family of matrices $\{A_i\}_{i=1}^{k}$ such that
\begin{align*}
    D_x f(x) \in \mathrm{conv}\{A_1,\ldots,A_k\},\quad  \mbox{ for every }\;  x\in \real^n.
\end{align*}
We assume that there exists a safe region in the state space denoted by $\mathcal{X}\subset \real^n$. The goal is to design a state feedback controller $u=Fx$ for the control system~\eqref{eq:safety} such that the closed-loop system avoids the unsafe region in the state space, for any bounded disturbances $w(t)\in [\underline{w},\overline{w}]$. We also assume that these exists a polytope $\mathcal{P}\subseteq \real^n$: 
\begin{align*}
   \mathcal{P}=\setdef{x\in \real^n}{\underline{h}\le Hx \le \overline{h}},
\end{align*}
where $H\in\real^{n\times m}$ is a full column rank matrix and $\underline{h}\le \overline{h}\in \real^m$ are such that the safe set $\mathcal{X}$ can be under-approximated by $\mathcal{P}$, i.e., we have $\mathcal{P}\subseteq \mathcal{X}$. Using the polytope $\mathcal{P}$, one can define a polyhedral cone $K\subseteq\real^n$ with the following $H$-rep: 
\begin{align}\label{eq:hrep-app}
    K =\setdef{x\in \real^n}{ Hx \ge \vect{0}_m}.
\end{align}
Let $V\in \real^{n\times q}$ be a generating matrix for the cone $K$, i.e., $K$ has a $V$-rep given by $K=\setdef{Vx}{x\ge \vect{0}_q}$. Indeed, using the cone $K$, the polytope $\mathcal{P}$ can be described by the interval $[\underline{\eta},\overline{\eta}]_K$, where $\underline{\eta},\overline{\eta}\in \real^n$ is such that $H\underline{\eta}=\underline{h}$ and $H\overline{\eta}=\overline{h}$. We introduce the following linear programming feasibility problem with unknown parameters $F$ and $\alpha$:
\begin{align}\label{eq:LP}
    &H(A_i+BF+\alpha I_n)V \ge \vect{0}_{m\times q} ,\quad\forall i \in \{1,\ldots,k\}\nonumber\\
    &H(f(\overline{\eta}) + BF\overline{\eta}) + (HC)^+ \overline{w} + (HC)^{-}\underline{w} \le \vect{0}_m, \nonumber\\
    &H(f(\underline{\eta}) + BF\underline{\eta}) + (HC)^+\underline{w} + (HC)^{-}\overline{w} \ge \vect{0}_m.
    \end{align}

\begin{theorem}[Control design via linear programming]\label{thm:LP}
Consider dynamical system~\eqref{eq:safety} with the polyhedral cone $K\in \real^n$ defined in~\eqref{eq:hrep-app}. Suppose that linear programming~\eqref{eq:LP} is feasible with a solution  $(F^*,\alpha^*)$. By choosing the state feedback controller $u=F^*x$, we obtain the closed-loop system
\begin{align}\label{eq:closedloop}
    \dot{x} = f(x)+BF^*x + Cw.
\end{align}
Then the polytope $\mathcal{P}$ is a forward invariant set for the system~\eqref{eq:closedloop} for any disturbance $t\mapsto w(t)$  such that $w(t)\in [\underline{w},\overline{w}]$ for all $t\in \real_{\ge 0}$.
\end{theorem}
  \smallskip
\begin{proof}
Consider the closed-loop system~\eqref{eq:closedloop}. First, by the linear programming~\eqref{eq:LP}, $F^*\in \real^{n\times p}$ satisfies $H(A_i+BF^*+\alpha^*I_n)V\ge \vect{0}_{m\times q}$, for every $i\in \{1,\ldots,k\}$.  This implies that $H(D_xf(x)+BF^* + \alpha^*I_n)V\ge \vect{0}_{m\times q}$, for every $x\in \real^n$. Therefore, by Theorem~\ref{thm:monotone}\ref{p2:mon}, the closed-loop system~\eqref{eq:closedloop} is $K$-monotone. Thus, for every disturbance $t\mapsto w(t)$ with $w(t)\in [\underline{w},\overline{w}]$,
\begin{multline*}
    H(f(\overline{\eta})+BF^*\overline{\eta}+ Cw(t)) \\ \le H(f(\overline{\eta}) + BF^*\overline{\eta}) + (HC)^+ \overline{w} + (HC)^{-}\underline{w} \le \vect{0}_m
\end{multline*}
where the first inequality holds because $w(t)\in [\underline{w},\overline{w}]$ and the second inequality holds by the linear programming~\eqref{eq:LP}. Therefore, by~\cite[Proposition 2.1]{HLS:95}, the trajectory of the closed-loop system starting from $\overline{\eta}$ is non-increasing with respect to the preorder $\preceq_K$. This means that $\setdef{x\in \real^n}{x\preceq_K \overline{\eta}}$ is invariant for the closed-loop system~\eqref{eq:closedloop}. Similarly, one can use the constraint $H(f(\underline{\eta})+ BF^*\underline{\eta} + Cw(t))\ge \vect{0}_m$ to show that $\setdef{x\in \real^n}{\underline{\eta}\preceq_K x}$ is an invariant set for the closed-loop system~\eqref{eq:closedloop}. As a result, $\setdef{x\in \real^n}{\underline{\eta} \preceq_K x\preceq_K \overline{\eta}} = [\underline{\eta},\overline{\eta}]_K=\mathcal{P}$ is an invariant set for the closed-loop system~\eqref{eq:closedloop}. 
\end{proof}

\begin{example}[Feedback Controller for Inverted Pendulum]\label{ex:pendulum}
Consider the inverted pendulum with the following dynamics:
\begin{align}\label{eq:ip}
    \dot{x}_1&=x_2,\nonumber\\
    \dot{x}_2&=\frac{g}{\ell} \sin(x_1) + u + w,
\end{align}
where $x_1$ and $x_2$ are the angular position and the angular velocity of the pendulum, $u\in \real$ is the control input, and $w\in \real$ is the disturbance. In this example $g$ is the gravitational constant and $\ell$ is the length of the pendulum.
\begin{figure}
 \begin{center}
		\includegraphics[width =0.6\linewidth,clip]{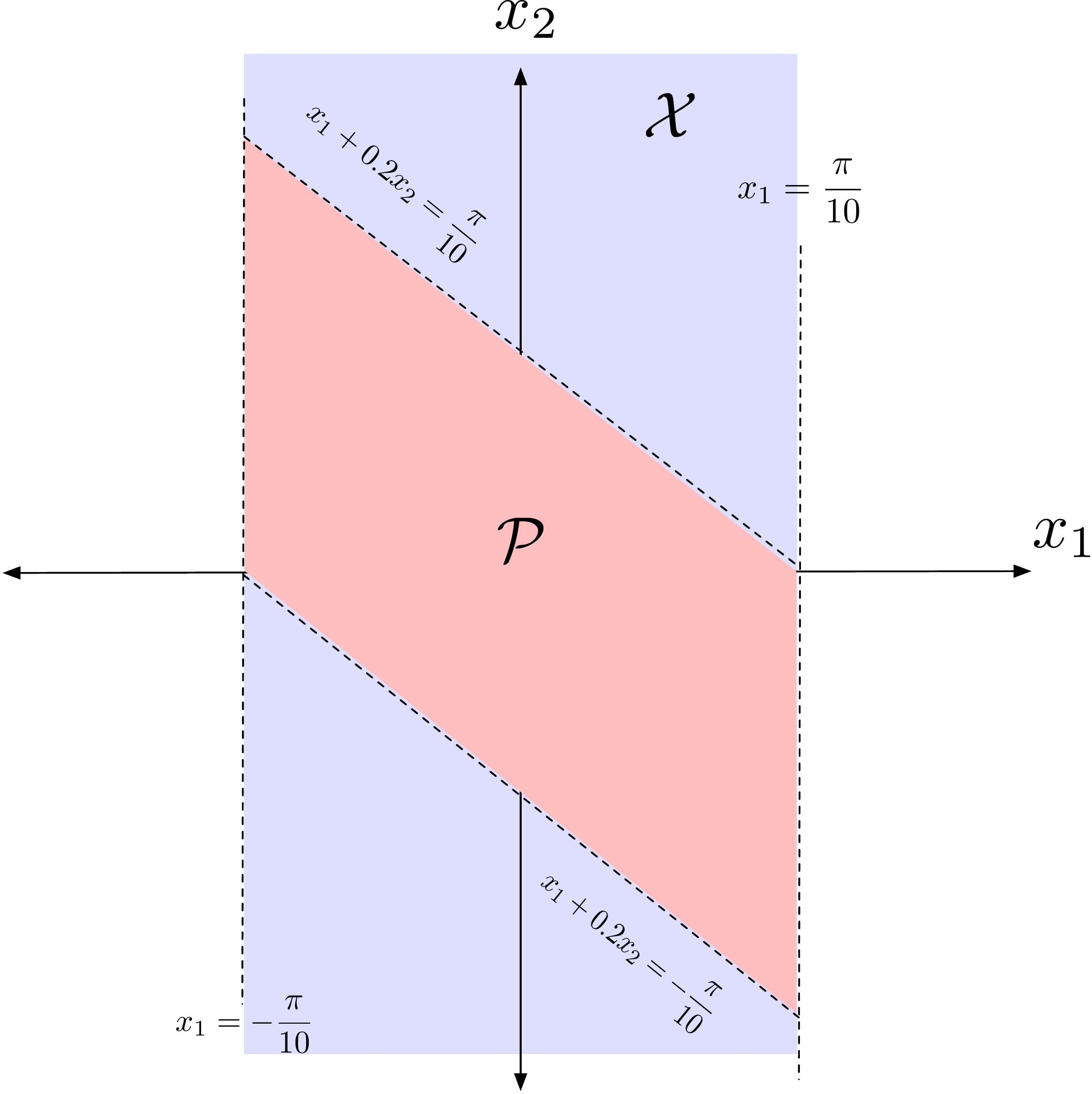}
	\end{center}
	\caption{The safe set $\mathcal{X}$ (blue) and its under-approximation by the polytope $\mathcal{P}$ (red).}
	\label{fig:pendulum}
	\vspace{-0.3cm}
\end{figure}
We assume that $\frac{g}{\ell}=1$ and the disturbance is a time-varying unknown signal with $w(t)\in [-0.2,0.2]$, for every $t\in \real_{\ge 0}$. The safe set is given by $\mathcal{X}=\setdef{(x_1,x_2)\in \real^2}{\frac{-\pi}{10}\le x_1\le \frac{\pi}{10}}$ and is shown in black in Figure~\ref{fig:pendulum}. First, note that $x_1=x_2=0$ is an unstable equilibrium point for the inverted pendulum~\eqref{eq:ip} without any input or disturbances. Second, for the inverted pendulum~\eqref{eq:ip}, there exists no controller that can make a rectangular neighborhood of $\vect{0}_2$ forward invariant. This can be proved as follows: consider a rectangular neighborhood of $\vect{0}_2$ and assume that $x_1=\lambda>0$ is an edge of this rectangle. Since the neighborhood contains $\vect{0}_2$, there are points on this edge such that $x_2>0$. This implies that on this edge, we have $\dot{x}_1=x_2>0$. Thus, this rectangular set cannot be forward invariant for the system~\eqref{eq:ip}. Now, 
we consider the under-approximation of the safe set $\mathcal{X}$ by the polytope $\mathcal{P}$ as shown in red in Figure~\eqref{fig:pendulum} and described by
\begin{align*}
    \mathcal{P} = \setdef{(x_1,x_2)^{\top}\in \real^2}{-\tfrac{\pi}{10}\left[\begin{smallmatrix}1\\ 1\end{smallmatrix}\right] \le \left[\begin{smallmatrix}1 & 0.2\\ 1 & 0\end{smallmatrix}\right]\left[\begin{smallmatrix}x_1\\ x_2\end{smallmatrix}\right]\le \tfrac{\pi}{10}\left[\begin{smallmatrix}1\\ 1\end{smallmatrix}\right] }.
\end{align*}
One can define the cone $K=\setdef{(x_1,x_2)^{\top}\in \real^2}{\left[\begin{smallmatrix}1 & 0.2\\ 1 & 0\end{smallmatrix}\right]\left[\begin{smallmatrix}x_1\\ x_2\end{smallmatrix}\right]\ge \vect{0}_2}$ and check that it has the following $V$-rep: $K=\setdef{\left[\begin{smallmatrix}0 & 1\\ 5 & -5\end{smallmatrix}\right]\left[\begin{smallmatrix}x_1\\ x_2\end{smallmatrix}\right]}{(x_1,x_2)^{\top}\in \real^2_{\ge 0}}$. Additionally, one can solve $\left[\begin{smallmatrix}1 & 0.2\\ 1 & 0\end{smallmatrix}\right]\overline{\eta}=-\left[\begin{smallmatrix}1 & 0.2\\ 1 & 0\end{smallmatrix}\right]\underline{\eta}=\frac{\pi}{10}\left[\begin{smallmatrix}1\\ 1\end{smallmatrix}\right]$ to obtain $\overline{\eta}=-\underline{\eta}=[\frac{\pi}{10}, \;\; 0]^{\top}$. Moreover, the Jacobian of the inverted pendulum at $(x_1,x_2)^{\top}$ is given by $\left[\begin{smallmatrix} 0 & 1 \\ \frac{g}{\ell}\cos(x_1) & 0\end{smallmatrix}\right]$ and, for every $(x_1,x_2)^{\top}\in \mathcal{P}$,
\begin{align*}
    \left[\begin{smallmatrix}  0 & 1 \\ \frac{g}{\ell}\cos(x_1) & 0\end{smallmatrix}\right] \in \mathrm{conv}\left\{A_1:=\left[\begin{smallmatrix}  0 & 1 \\ \frac{g}{\ell} & 0\end{smallmatrix}\right],A_2:=\left[\begin{smallmatrix}  0 & 1 \\ \frac{g}{\ell}\cos(\frac{\pi}{10}) & 0\end{smallmatrix}\right]\right\}.
\end{align*}
The optimal solution of the linear program~\eqref{eq:LP} is $(F^*,\alpha^*) = (\left[\begin{smallmatrix}1.6203 \\ 
    5.1338\end{smallmatrix}\right], 20)$. Thus, by applying the feedback controller $u = -1.6203 x_1 - 5.1338 x_2$ and using Theorem~\ref{thm:LP}, the polytope $\mathcal{P}$ is  forward invariant for any disturbance in the interval $[-0.2,0.2]$. 
\end{example}

\section{Conclusions}

We characterize monotonicity of a control system with respect to a polyhedral cone using the half-space representation and the vertex representation of the cone. We use the notion of gauge norm as a key element for connecting contraction theory with monotone theory on polyhedral cones. We provide necessary and sufficient conditions for contractivity of monotone systems with respect to the gauge norms.

\section*{Appendices}
\addcontentsline{toc}{section}{Appendices}
\renewcommand{\thesubsection}{\Alph{subsection}}

\subsection{Proof of Lemma~\ref{thm:pos-mon}}\label{app:B}
             Regarding~\ref{p1:positive} $\implies$~\ref{p2:monotone}, since $I_n + hA$ is $K$-positive, for every $\phi\in K^*$ and every $x\in K$ such that $\langle\phi,x\rangle=0$, we have $(I_n + h A)x\in K$. This implies that $0\le \langle\phi,(I_n + h A)x\rangle  = h \langle\phi,Ax\rangle$. 
             Since $h>0$, we have $\langle\phi,Ax\rangle\ge 0$ and thus $A$ is $K$-Metzler. Regarding~\ref{p2:monotone} $\implies$~\ref{p3:positive-diag}, suppose that $K$ is a polyhedral cone with generating linear functionals $\{\phi_i\}_{i=1}^{m}$. {\color{black}Suppose the $V$-rep of $K$ is given by $K=\setdef{Vx}{x\ge \vect{0}_q}$ where $v_j\in \real^n$ is the $j$th column of the generating matrix $V\in \real^{n\times q}$, for $j\in\{1,\ldots,q\}$. For every $i\in \{1,\ldots,m\}$, we set $\mathcal{S}_{i} = \setdef{j\in \{1,\ldots,q\}}{\langle\phi_i,v_j\rangle \ne 0}$. First, assume that $\mathcal{S}_i\ne \emptyset$ and define the constant $\alpha_i^* = \max_{j\in \mathcal{S}_i}\frac{-\langle\phi_i,Av_j\rangle}{\langle\phi_i,v_j\rangle}$.
             Then, for every $j\in \mathcal{S}_i$ and every $\alpha\ge \alpha_i^*$,
             \begin{align}\label{eq:1}
                 \langle\phi_i,(A+\alpha I_n)v_j\rangle =  \langle\phi_i,A v_j\rangle +  \alpha \langle\phi_i,v_j\rangle \ge 0,
             \end{align}
             where the second inequality holds by definition of $\alpha_i^*$ and the fact that $\alpha\ge \alpha_i^*$.  
             Moreover, for every $j\not\in \mathcal{S}_i$ and every $\alpha\ge \alpha_i^*$, 
             \begin{align}\label{eq:2}
                 \langle\phi_i,(A+\alpha I_n)v_j\rangle = \langle\phi_i,A v_j\rangle  \ge 0,
             \end{align}
             where the first equality holds because $\langle\phi_i,v_j\rangle = 0$, for $j\not\in\mathcal{S}_i$, and the second inequality holds because $A$ is $K$-Metzler. Now, assume $\mathcal{S}_i=\emptyset$. In this case, for every $j\in \{1,\ldots,q\}$ and every $\alpha\in \real$,
             \begin{align*}
                 \langle\phi_i,(A+\alpha I_n)v_j\rangle = \langle\phi_i,A v_j\rangle  \ge 0,
             \end{align*}
             where the first equality holds because $\langle\phi_i,v_j\rangle = 0$, for every $j\in \{1,\ldots,q\}$ and the second inequality holds because $A$ is $K$-Metzler. Since $x\in K$, there exist constants $c_1,\ldots,c_q\ge 0$ such that $x=\sum_{j=1}^{q}c_j v_j$. This implies that, for every $\alpha\ge \alpha_i^*$, we have
             \begin{align*}
                 \langle\phi_i,(A+\alpha I_n)x\rangle & = \sum\nolimits_{j\in \mathcal{S}_i} c_j \langle\phi_i,(A+\alpha I_n)v_j\rangle \\ & + \sum\nolimits_{j\not\in \mathcal{S}_i}c_j \langle\phi_i,(A+\alpha I_n)v_j\rangle \ge 0,
             \end{align*} 
             where the last inequality holds using inequalities~\eqref{eq:1} and~\eqref{eq:2}. Therefore, for every $i\in \{1,\ldots,m\}$, every $x\in K$, and every $\alpha\ge \max_{i\in \{1,\ldots,m\}}\alpha^*_i$, we have 
            \begin{align*}
              \langle\phi_i,(A+\alpha I_n)x\rangle \ge 0.
            \end{align*}
            As a result, for every $\alpha\ge \max_{i\in \{1,\ldots,m\}}\alpha_i^*$, we have $(A+\alpha I_n)K\subseteq K$. This means that, for every $\alpha\ge \max_{i\in \{1,\ldots,m\}}\alpha_i^*$, the matrix $A+\alpha I_n$ is $K$-positive. Regarding ~\ref{p3:positive-diag}$\implies$~\ref{p1:positive}, we first use the identity $A+\alpha I_n =  \alpha(I_n + \frac{1}{\alpha}A)$. The matrix $A + \alpha I_n$ is $K$-positive, for every $\alpha\ge \alpha^*$. By setting $h^*=\frac{1}{\alpha^*}$, this implies that $I_n + hA = h(A + \frac{1}{h}I_n)$ is $K$-positive, for every $0\le h\le h^*=\frac{1}{\alpha^*}$.}
          
\subsection{Proof of Theorem~\ref{thm:contraction}}\label{ap:C}
             Regarding~\ref{p1} $\implies$~\ref{p2}, let $\epsilon>0$ and note that by
             definition of the matrix measure, we have
             \begin{align}\label{eq:similar}
              &\lim_{h\to 0^+}\frac{\|(I_n+h A)\mathbf{e}\|_{\mathbf{e},K}-1}{h}
                \le \lim_{h\to 0^+}\frac{\sup_{v\ne
              0}\frac{\|(I_n+h A)v\|_{\mathbf{e},K}}{\|v\|}-1}{h} \nonumber\\ &
              \le \lim_{h\to 0^+}\frac{\|I_n+h A\|_{\mathbf{e},K}-1}{h} =
              \mu_{\mathbf{e},K}(A) \le c < c + \epsilon.
             \end{align}
             Note that the function $h\mapsto
             \frac{\|(I_n+h A)\mathbf{e}\|_{\mathbf{e},K}-1}{h}$ is a
             {\color{black}weakly increasing function on $(0,\infty)$.} 
             This implies that there exists $h^*> 0$ such that
             $\frac{\|(I_n+h A)\mathbf{e}\|_{\mathbf{e},K}-1}{h} <
             c+\epsilon$ for every $h\in [0,h^*]$. As a result, we get
             $\|(I_n+h A)\mathbf{e}\|_{\mathbf{e},K} \le 1+(c+\epsilon)h$, for every $h\in [0,h^*]$. Using the definition of the gauge norm, for every $h\in [0,h^*]$, 
               \begin{align*}
                 -(1+(c+\epsilon) h)\mathbf{e} \preceq_K (I_n+h A)\mathbf{e}
                 \preceq_K (1+(c+\epsilon) h)\mathbf{e}.
               \end{align*}
               This means that $A\mathbf{e} \preceq_K (c+\epsilon) \mathbf{e}$. Since the LHS of this inequality is independent of $\epsilon$, we can deduce that $A\mathbf{e} \preceq_K c \mathbf{e}$. Regarding~\ref{p2} $\implies$~\ref{p1}, suppose that
               $v\in \real^n$ is such that
               $\|v\|_{\mathbf{e},K}=\lambda$. This means that
               $\lambda$ is the smallest positive number such that $-\lambda\mathbf{e} \preceq_K v \preceq_K \lambda \mathbf{e}$.
               Note that by Theorem~\ref{thm:pos-mon}, there exists $h^*>0$ such
               that $I_n+h A$ is $K$-positive for every $h\in
               [0,h^*]$. As a result, for every $h\in [0,h^*]$, we
               have $ \vect{0}_n\preceq_K (I_n + h A)(\lambda \mathbf{e} -
                 v)$ and this implies that $(I_n+h A)v \preceq_K \lambda (I_n + h
               A)\mathbf{e}$, for every $h\in [0,h^*]$. Similarly, one can show that $-\lambda (I_n + h
               A)\mathbf{e} \preceq_K (I_n+h A)v$, for every $h\in
               [0,h^*]$.  Thus, for every $h\in [0,h^*]$,
               \begin{align*}
                 -\lambda(1+c h) \mathbf{e} &\preceq_K -\lambda (I_n +
                 h A)) \mathbf{e}  \preceq_K (I_n +h A)v \\ &\preceq_K \lambda (I_n +
                 h A) \mathbf{e} \preceq_K \lambda(1+c h).
               \end{align*}
               This means that, for every $h\in [0,h^*]$, we have $|(I_n+h A)v\|_{\mathbf{e},K} \le
                 \|v\|_{\mathbf{e},K}(1+c h)$.
               Using the definition of the matrix measure, we get
               $\mu_{\mathbf{e},K}(A)\le c$. Regarding part~\ref{p2} $\iff$~\ref{p3}, the result follows from
               Lemma~\ref{thm:inequality}. Regarding~\ref{p1:dual} $\implies$~\ref{p2:dual}, consider $v\in K$ such that $\|v\|^{\mathrm{d}}_{\mathbf{e}^*,K^*}=1$. For every $\eta\in [-\mathbf{e}^*,\mathbf{e}^*]$, we have $ \mathbf{e}^*-\eta \in K^*$ and $\mathbf{e}^*+\eta \in K^*$. This implies that $-\langle\mathbf{e}^*,v\rangle \le \langle\eta,v\rangle \le \langle\mathbf{e}^*,v\rangle$. 
               Now, using the definition of the dual norm, we get
               \begin{align*}
                  \langle\mathbf{e}^*,v\rangle = \max \setdef{|\langle \eta, v\rangle|}{\eta\in [-\mathbf{e}^*,\mathbf{e}^*]_{K^*}}= \|v\|^{\mathrm{d}}_{\mathbf{e}^*,K^*}=1.
               \end{align*}
               Using the fact that $\|v\|^{\mathrm{d}}_{\mathbf{e}^*,K^*}=1$, a similar argument as in~\eqref{eq:similar} will lead to the following inequality:
               \begin{align*}
               \lim_{h\to 0^+}&\frac{\|(I_n+h A)v\|^{\mathrm{d}}_{\mathbf{e}^*,K^*}-1}{h}
               \le 
               \mu^{\mathrm{d}}_{\mathbf{e}^*,K^*}(A) \le c.
             \end{align*}
             Therefore, we have $\|(I_n+h A)v\|^{\mathrm{d}}_{\mathbf{e}^*,K^*} \le 1+ch$ and thus, by definition of the dual gauge norm,
             \begin{align*}
                 -1-ch \le \langle \mathbf{e}^*, (I_n+h A)v\rangle\le  1+ch. 
             \end{align*}
             As a result, for every $v\in K$, such that $\|v\|^{\mathrm{d}}_{\mathbf{e}^*,K^*}=1$,
             \begin{align*}
                 -1-ch \le \langle (I_n+h A^{\top})\mathbf{e}^*, v\rangle\le  1+ch. 
             \end{align*}
             Using the fact that $\langle \mathbf{e}^*,v\rangle=1$, we get
             $-c \langle \mathbf{e}^*,v\rangle \le \langle A^{\top}\mathbf{e}^*, v\rangle\le c \langle \mathbf{e}^*,v\rangle$. Note that the above inequalities hold for every $v\in K$ satisfying $\|v\|^{\mathrm{d}}_{\mathbf{e}^*,K^*}=1$. Therefore, by definition of the preorder $\preceq_{K^*}$, we get $ A^{\top}\mathbf{e}^* \preceq_{K^*} c\mathbf{e}^*$. Regarding~\ref{p2:dual} $\implies$~\ref{p1:dual}, suppose that $\phi$ is such that $-\mathbf{e}^*\preceq_{K^*}\phi\preceq_{K^*}\mathbf{e}^*$. Since $A$ is $K$-Metzler, by Lemma~\ref{thm:pos-mon}, there exists $h^*>0$ such
               that $I_n+h A$ is $K$-positive for every $h\in
               [0,h^*]$. Therefore, using~\cite[Theorem 2.24]{AB-RJP:94}, the operator $I_n+h A^{\top}$ is $K^*$-positive, for every $h\in
               [0,h^*]$. As a result, for every $h\in [0,h^*]$, we
               have $ \vect{0}_n\preceq_{K^*} (I_n + h A^{\top})(\mathbf{e}^* -
                 \phi)$ and this implies that $(I_n+h A^{\top})\phi \preceq_{K^*}  (I_n + h
               A^{\top})\mathbf{e}^*$, for every $h\in [0,h^*]$. Similarly, one can show that $- (I_n + h
               A^{
               \top})\mathbf{e}^* \preceq_{K^*} (I_n+h A^{\top})\phi$, for every $h\in
               [0,h^*]$.  This implies that, for every $h\in [0,h^*]$,
               \begin{align*}
                 -(1+c h) \mathbf{e}^* &\preceq_{K^*} - (I_n +
                 h A^{\top}) \mathbf{e}^*  \preceq_{K^*} (I_n +h A^{\top})\phi \\ &\preceq_{K^*}  (I_n +
                 h A^{\top}) \mathbf{e}^* \preceq_{K^*} (1+c h)\mathbf{e}^*.
               \end{align*}
               Therefore, for every $-\mathbf{e}^*\preceq_{K^*} \phi \preceq_{K^*} \mathbf{e}^*$, we get
               \begin{align*}
                   \|(I_n+h A)v\|&_{\mathbf{e}^*,K^*} = \max |\langle \phi, (I_n+h A)v\rangle | \\ & = \max |\langle (I_n+h A^{\top})\phi, v\rangle| \le (1+ch) \|v\|^{\mathrm{d}}_{\mathbf{e}^*,K^*}
               \end{align*}
               where the last equality holds by definition of the dual gauge norm. Using the definition of the matrix measure, we get $\mu^{\mathrm{d}}_{\mathbf{e}^*,K^*}(A)\le c$. Regarding~\ref{p2:dual} $\iff$~\ref{p3:dual}, the result follows from Lemma~\ref{thm:inequality}.

\bibliographystyle{IEEEtran}
\bibliography{SJ.bib}
\end{document}